\documentclass[11pt]{amsart}
\usepackage{latexsym,amsfonts,amssymb,amsmath,amsthm, verbatim}
\usepackage[all]{xy}

\newcommand\Ell[2]{\ell(#1 | #2)}
\newcommand{\pos}{\mathrm{pos}}
\newcommand{\qcomm}[3]{q^{{#1}\langle\!\langle #2, #3\rangle\!\rangle}}
\newcommand\qdet[1]{[\![#1]\!]}
\newcommand{\detq}{\mathrm{det}_q}

\newcommand{\newatop}[2]{\genfrac{}{}{0pt}{}{#1}{#2}}

\newcommand{\llangle}{\langle\!\langle}
\newcommand{\rrangle}{\rangle\!\rangle}
\newcommand{\Sym}{\mathfrak S}
\newcommand{\mM}{\mathsf{mM}}
\newcommand{\dcup}{\,\dot{\cup}\,}
\newcommand{\surr}{\mathrel{\hbox{$\curvearrowright$\hspace{-1.45ex}\rule[0.1ex]{.22ex}{.22ex}\hspace{1.23ex}}}}

\theoremstyle{plain}		
	\newtheorem{mytheo}{Theorem}
	\newtheorem{myprop}[mytheo]{Proposition}
	\newtheorem{mycoro}[mytheo]{Corollary}

\theoremstyle{definition}	
	\newtheorem{mydefi}{Definition}

\theoremstyle{remark}	
	\newtheorem{mylemm}{Lemma}
	\newtheorem{myrema}{Remark}
	\newtheorem*{myexam}{Example}

\begin{document}

\title{Quasideterminants and $q$-commuting minors}
\author{Aaron Lauve}
\date{February 09, 2006}

\begin{abstract}
We present two new proofs of the the important $q$-commuting property holding among certain pairs of quantum minors of an $n\times n$ $q$-generic matrix. The first uses elementary quasideterminantal arithmetic; the second involves paths in an edge-weighted directed graph.
\end{abstract}

\maketitle

\def\WeaklySeparated{In the literature, sets $J$ and $I$ sharing this relationship are called ``weakly separated.'' We avoid this terminology because it does not indicate who separates whom.}
\section{Introduction \& Main Theorem}\label{sec:intro}
This paper arose from an attempt to understand the ``quantum shape algebra'' of Taft and Towber \cite{TafTow:1}, which we call the \emph{flag algebra} $\mathcal{F}\ell_q(n)$ here. One goal was to find quasideterminantal justifications for the relations chosen for $\mathcal{F}\ell_q(n)$. A second goal was to find some hidden relations, within $\mathcal{F}\ell_q(n)$, known to hold in an isomorphic image. To more quickly reach a statement of the theorem, we save further remarks on the goals for later.

\begin{mydefi}\label{def:surrounds} Given two subsets $I,J\subseteq[n]$, we say $J$ \emph{surrounds}\footnote{\WeaklySeparated} $I$, written $J\surr I$, if (i) $|J|\leq |I|$, and (ii) there exist disjoint subsets $\emptyset\subseteq J',J''\subseteq J$ such that:
\begin{enumerate}
\item[a.] $J\setminus I = J' \dcup J''$,
\item[b.] $j'<i$ for all $j'\in J'$ and $i\in I\setminus J$,
\item[c.] $i<j''$ for all $i\in I\setminus J$ and $j''\in J'$,
\end{enumerate}
In this case, we put $\llangle J,I\rrangle = |J''|-|J'|$.
\end{mydefi}

Given an $n\times n$ \emph{$q$-generic matrix} $X$ and a subset $I\subseteq[n]$ with $|I|=d$, we write $\qdet{I}$ for the quantum minor built from $X$ by taking row-set $I$ and column-set $[d]$. 

\begin{mytheo}[$q$-Commuting Minors]\label{thm:qComm} If the subsets $I,J\subseteq[n]$ satisfy $J\surr I$, the quantum minors $\qdet{J}$ and $\qdet{I}$ \emph{$q$-commute}. Specifically,
\begin{equation}\label{eq:qComm}
\qdet{J}\qdet{I} = q^{\llangle J,I\rrangle} \qdet{I}\qdet{J}\,.
\end{equation}
\end{mytheo}

An earlier proof of this theorem may be found in \cite{KroLec:1}, while Leclerc and Zelevinsky \cite{LecZel:1} show that if $\qdet{J} \qdet{I} = q^{\alpha} \qdet{I} \qdet{J}$ for some $\alpha\in\mathbb{Z}$, then $J\surr I$. We give two new proofs in the sequel. The first proof ($\mathcal{Q}$) uses simple arithmetic involving quasideterminants; the second ($\mathcal{G}$) involves counting weighted paths on a directed graph.

\subsection{Useful notation}\label{sec:notation}
The reader has already encountered our notation $[n]$ for the set $\{1,2,\ldots, n\}$; let $\binom{[n]}{d}$ denote the set of all subsets of $[n]$ of size $d$. Given a set $I=\{i_1<i_2<\cdots i_d\}\in\binom{[n]}{d}$ and any $I'\subseteq I$, we write $I^{I'}$ for the subset built from $I$ by \emph{deleting} $I'$ (i.e. $I\setminus I'$) and $I_{I'}$ for the complement (i.e. a fancy way of saying \emph{keep} $I'$). In case $\Lambda = \{\lambda_1< \lambda_2 \cdots <\lambda_r\} \in\binom{[d]}{r}$, we write $I_{(\Lambda)}$ for the subset $\{i_{\lambda_1}, i_{\lambda_2}, \ldots, i_{\lambda_r}\}$ and $I^{(\Lambda)}$ for the complement.

Suppose instead $I\in [n]^d$, the set of all $d$-tuples chosen from $[n]$. In this case, the notations $I_{I'}$ and $I^{I'}$ are not well-defined (as the entries of $I'$ may occur in more than one place within $I$) but the notations $I_{(\Lambda)}$ and $I^{(\Lambda)}$ will be useful in the sequel. If $I,J$ are two sets or tuples of sizes $d,e$ respectively, we define $A|B$ to be the $(d+e)$-tuple $(i_1,\ldots, i_d,j_1,\ldots, j_e)$. Let $[n]^d_{\ast}\subseteq[n]^d$ denote those $d$-tuples with distinct entries. For $I\in[n]^d_{\ast}$, we define the \emph{length} of $I$ to be $\ell(I) = \#\mathrm{inv}(I) = \#\big\{(j,k) : j<k\hbox{ and }i_j>i_k \big\}$. Fix $i\in[n]$ and $I = i_1,i_2,\ldots, i_d$ (viewed either as a set or a $d$-tuple without repetition); if there is a $1\leq k\leq d$ with $i_k = i$, then $k$ is the \emph{position} of $i$ and we write $\mathrm{pos}_I(i) = k$.

We extend our delete/keep notation to matrices. Let $A$ be an $n\times n$ matrix whose rows and columns are indexed by $R$ and $C$, respectively. For any $R'\subseteq R$ and $C'\subseteq C$, we let $A^{R',C'}$ denote the submatrix built from $A$ by deleting row-indices $R'$ and column-indices $C'$. Let $A_{R',C'}$ be the complementary submatrix. In case $R'=\{r\}$ and $C'=\{c\}$, we may abuse notation and write, e.g., $A^{rc}$. We will also need a means to construct matrices from $A$ whose rows (columns) are repeated or are not in their natural order. If $I\in R^{\,d}$ and $J\in C^{\,e}$, let $A_{I,J}$ denote the obvious new matrix built from $A$. 

\section{Preliminaries for $\mathcal Q$-Proof}\label{sec:prelims-Qproof}
\subsection{Quasideterminants}\label{sec:quasidets} 
The quasideterminant \cite{GelRet:1,GGRW:1} was introduced by Gelfand and Retakh as a replacement for the determinant over noncommutative rings $\mathcal R$. Given an $n\times n$ matrix $A=(a_{ij})$ over $\mathcal R$, the quasideterminant $|A|_{ij}$ (there is one for each position $(i,j)$ in the matrix) is not polynomial in the entries $a_{ij}$ but rather a rational expression, as we will soon see. Consequently, quasideterminants are not always defined. Below is a sufficient condition (cf. loc. cit. for more details).

\begin{mydefi} Given $A$ and $\mathcal R$ as above, if $A^{ij}$ is invertible over $\mathcal R$, then the $(ij)$-\emph{quasideterminant} is defined and given by
\[
|A|_{ij} = a_{ij} - \rho_i \cdot (A^{ij})^{-1} \cdot \chi_j\,,
\]
where $\rho_i$ is the $i$-th row of $A$ with column $j$ deleted and $\chi_j$ is the $j$-th column of $A$ with row $i$ deleted.
\end{mydefi}

\begin{myrema}\label{rem:quasidet-is-inv}
One deduces that $|A|_{ij}^{-1} = (A^{-1})_{ji}$ when both sides are defined. 
\end{myrema}

Details on this remark and the following three theorems may be found in \cite{GGRW:1,KroLec:1,Lau:1}. Note that the phrase `when defined' is implicit throughout.

\begin{mytheo}[Homological Relations]\label{thm:homol-rels} Let $A$ be a square matrix and let $i\neq j$ $(k\neq l)$ be two row (column) indices. We have
\[
-|A^{jk}|_{il}^{-1}\cdot|A|_{ik} =
|A^{ik}|_{jl}^{-1}\cdot|A|_{jk}.
\]
\end{mytheo}

\begin{mytheo}[Muir's Law of Extensible Minors]\label{thm:muir's-law} Let $A=A_{R,C}$ be a square matrix with row (column) indices $R$ ($C$). Fix $R_0\subsetneq R$ and $C_0\subsetneq C$. Say an algebraic, rational expression $\mathcal I=\mathcal I(A,R_0,C_0)$ involving the quasi-minors $\left\{ |A_{R',C'}|_{rc} : r\in R'\subseteq R_0\right.$, $\left.c\in C'\subseteq C_0\right\}$ is an \emph{identity} if the equation $\mathcal I = 0$ is valid. For any $L\subseteq R\setminus R_0$ and $M\subseteq C\setminus C_0$, the expression $\mathcal I'$ built from $\mathcal I$ by extending all minors $|A_{R',C'}|_{rc}$ to $|A_{L\cup R',M\cup C'}|_{rc}$ is also an identity.
\end{mytheo}

\begin{mydefi} Let $B$ be an $n\times d$ matrix. For any $i,j,k\in[n]$ and $M\subseteq[n]\setminus \{i\}$ ($|M|=d-1$), define $r_{ji}^M=r_{ji}^M(B):= |B_{(j|M), [d]}|_{jk}|B_{(i|M), [d]}|_{ik}^{-1}$. Gelfand and Retakh \cite{GelRet:3} show this ratio is independent of $k$, and call it a \emph{right-quasi-Pl\"ucker} coordinate for $B$. 
\end{mydefi} 

\begin{myrema} In case $B$ is $n\times m$ for some $m>d$, we choose the first $d$ columns of $B$ to form the above ratio unless otherwise indicated.
\end{myrema}

\begin{mytheo}[Quasi-Pl\"ucker Relations]\label{thm:qplucker-rels} Fix an $n\times n$ matrix $A$, subsets $M,L\subseteq[n]$ with $|M|+1\leq |L|$, and $i\in[n]\setminus M$. We have
the \emph{quasi-Pl\"ucker relation} $(\mathcal P_{L,M,i})$
\[
1 = \sum_{j\in L} r_{ij}^{L\setminus j} r_{ji}^{M}\,.
\]
\end{mytheo}

\subsection{Quantum determinants}\label{sec:quantumdets} An $n\times n$ matrix $X=(x_{ab})$ is said to be $q$-generic if its entries satisfy the relations 
\begin{eqnarray*}
(\forall i,\, \forall k<l)\,\quad x_{il}x_{ik} &=& q x_{ik} x_{il} \\
(\forall i<j,\, \forall k)\quad x_{jk}x_{ik} &=& q x_{ik} x_{jk} \\
(\forall i<j,\, \forall k<l)\quad x_{jk}x_{il} &=& x_{il} x_{jk} \\
(\forall i<j,\, \forall k<l)\quad x_{jl}x_{ik} &=& x_{ik} x_{jl} + (q-q^{-1}) x_{il}x_{jk}\,.
\end{eqnarray*}
Notice that every submatrix of a $q$-generic matrix is again $q$-generic. 

Fix a field $\Bbbk$ of characteristic 0 and a distinquished invertible element $q\in\Bbbk$ not equal to a root of unity. Let $\mathrm{M}_q(n)$ be the $\Bbbk$-algebra with $n^2$ generators $x_{ab}$ subject to the relations making $X$ a $q$-generic matrix. It is known \cite{KelLenRig:1} that $\mathrm{M}_q(n)$ is a (left) Ore domain with (left) field of fractions $D_q(n)$.

\begin{mydefi} Given any $d\times d$ matrix $A$, define $\detq A$ by
\[
\detq A = \sum_{\sigma\in\Sym_d} (-q)^{-\ell(\sigma)} a_{\sigma(1),1} a_{\sigma(2),2} \cdots a_{\sigma(d),d}\,.
\]
\end{mydefi}
When $A=X_{R,C}$ is a submatrix of $X$, we have: (i) this quantity agrees with the analogous quantity modeled after the column-permutation definition of the determinant, (ii) swapping two adjacent rows of $A$ introduces a $q^{-1}$, and (iii) allowing any row of $A$ to appear twice yields zero. Properties (ii) and (iii) allow us to uniquely define the determinant of $A=X_{I,C}$ for any $I\in[n]^d$ and $C\in\binom{[n]}{d}$. In case $C=\{1,2\ldots, d\}$, we introduce the shorthand notation $\detq A = \qdet{I}$. We will also need the case $C=s+[d]:=\{s+1,s+2, \ldots, s+d\}$ for some $s>0$, which we write as $\qdet{I;s}$. 

Properties (i)--(iii) give us the important

\begin{mytheo}[Quantum Determinantal Identities]\label{thm:qdet-identities} Let $A=X_{R,C}$ be a $d\times d$ submatrix of $X$. Then for all $i, j\in R$ and $k\in C$, we have: 
\begin{equation*}\label{eq:qLaplace}
\sum_{c\in C} A_{jc} \cdot \left\{(-q)^{\pos_I(i)-\pos_C(c)} \detq A^{{i}c}\right\} = \delta_{ij}\cdot \detq A
\end{equation*}
\begin{equation*}\label{eq:qdet-minors-commute}
\Big[\detq A \,\,,\,A_{ik}\Big] = 0\,.
\end{equation*}
\end{mytheo}

In particular every submatrix of $X$ is invertible in $D_q(n)$ and (after Remark \ref{rem:quasidet-is-inv}) we are free to use the preceding quasideterminantal formulas on matrices built from $X$. The important formula follows: for all $I\in[n]^d$
\begin{equation}\label{eq:quasidet-to-qdet}
|X_{I,\{s+1,\ldots, s+d\}}|_{i, s+d} = (-q)^{d-\pos_I(i)} \qdet{I;s} \cdot \qdet{I^{i};s}^{-1}\,,
\end{equation}
where the factors on the right commute. Theorems \ref{thm:homol-rels} and \ref{thm:qdet-identities} are combined with (\ref{eq:quasidet-to-qdet}) in \cite{Lau:1} to prove 

\begin{mytheo}\label{thm:weak-qComm} Given any $i,j\in[n]$, $\{j\}\surr \{i\}$. For any $M\subseteq[n]$, the quantum minors $\qdet{j| M}$ and $\qdet{i| M}$ $q$-commute according to equation (\ref{eq:qComm}).
\end{mytheo}

\section{$\mathcal Q$-Proof of Theorem}\label{sec:Qproof}
Our first proof of Theorem \ref{thm:qComm} proceeds by induction on $|J|$ and rests on two key lemmas.
\begin{mylemm}\label{thm:disjoint} If $I\subseteq [n]$ and $j\in [n]\setminus I$ satisfy $\{j\}\surr I$. Then $\qdet{j}\qdet{I} = \qcomm{}{j}{I}\qdet{I}\qdet{j}$.
\end{mylemm}
\begin{proof}
From $(\mathcal P_{I,\emptyset, j})$ and (\ref{eq:quasidet-to-qdet}) we have
\[
1= \sum_{i\in I} \qdet{j | I\setminus i}\qdet{i | I\setminus i}^{-1} \qdet{i}\qdet{j}^{-1} \,,
\]
or 
\begin{equation}\label{eq:toclear}
\qdet{j}= \sum_{i\in I} \qdet{j | I^i}\qdet{i | I^i}^{-1} \qdet{i} \,.
\end{equation}
Theorem \ref{thm:weak-qComm} tells us that $\qdet{j | I^i}$ and $\qdet{i | I^i}$ $q$-commute, so we may clear the denominator in (\ref{eq:toclear}) on the left and get
\begin{equation}\label{eq:clear-left}
\qdet{I}\qdet{j} = \sum_{i\in I}(-q)^{\Ell{i}{I^i}} \qcomm{-}{j}{I}\qdet{j|I^i}\qdet{i}\,.
\end{equation}
In the other direction, Theorem \ref{thm:qdet-identities} tells us that $\qdet{i|I^i}$ and $\qdet{i}$ commute; clearing (\ref{eq:toclear}) on the right yields
\begin{equation}\label{eq:clear-right}
\qdet{j}\qdet{I} = \sum_{i\in I}(-q)^{\Ell{i}{I^i}} \qdet{j|I^i}\qdet{i}\,.
\end{equation}
Compare (\ref{eq:clear-left}) and (\ref{eq:clear-right}) to conclude that $\qdet{j}$ and $\qdet{I}$ $q$-commute as desired.
\end{proof}

\begin{mylemm}\label{thm:overlap} Fix $J,I\subseteq [n]$ satisfying $J\surr I$. For all $M \subseteq [n]\setminus (I\cup J)$, one has $J\cup M \surr I\cup M$ and $\qdet{J\cup M}\qdet{I\cup M} = \qcomm{}{J}{I}\qdet{I\cup M}\qdet{J\cup M}$.
\end{mylemm}

\begin{proof} The first statement is clear from the definition of `surrounds.' The second statement is a consequence of Muir's Law. 

Let $J=\{j_1,\ldots, j_d\}$, $I=\{i_1,\ldots, i_e\}$, and $M=\{m_1,\ldots, m_s\}$. Because of the nature of the defining relations for $q$-generic matrices and the definition of quantum determinant, the expression $\qdet{J}\qdet{I} = \qcomm{}{J}{I}\qdet{I}\qdet{J}$ is equivalent to  $\qdet{J;s}\qdet{I;s} = \qcomm{}{J}{I}\qdet{I;s}\qdet{J;s}$, or even $\qdet{I;s}^{-1}\qdet{J;s} = \qcomm{}{J}{I}\qdet{J;s}\qdet{I;s}^{-1}$. 

Let us write the left-hand side of this last equation in terms of quasideterminants:
\begin{eqnarray*}
\qdet{I;s}^{-1}\qdet{J;s} &=& \Big(|X_{I,s+[e]}|_{i_e,s+e}^{-1} \,\qdet{I^{(e)};s}^{-1} \Big)\times \Big(\qdet{J^{(d)};s} \,|X_{J,s+[d]}|_{j_d,s+d}\Big) \\
&\vdots& \\
& = & |X_{I,s+[e]}|_{i_e,s+e}^{-1} \cdots |X_{(i_1,i_2),\{s+1,s+2\}}|_{i_2,s+2}^{-1} \,|X_{i_1,s+1}|_{i_1,s+1}^{-1} \,\times \\
&& |X_{j_1,s+1}|_{j_1,s+1}\, |X_{(j_1,j_2),\{s+1,s+2\}}|_{j_2,s+2} \cdots |X_{J,s+[d]}|_{j_d,s+d} .
\end{eqnarray*}
Do the same to the right-hand side and get an identity involving quasideterminants. Notice that the submatrix $X_{M,[s]}$ appears nowhere in that identity. Inserting this everywhere according to Muir's Law and multiplying and dividing by $\qdet{M}$ we get (for the left-hand side)
\begin{eqnarray*}
|X_{(I|M),[s+e]}|_{i_e,s+e}^{-1} \cdots |X_{(i_1|M),[s+1]}|_{i_1,s+1}^{-1} \,\qdet{M}^{-1}\,\times \\
\qquad \qdet{M}\, |X_{(j_1|M),[s+1]}|_{j_1,s+1} \cdots |X_{(J|M),[s+d]}|_{j_d,s+d} .
\end{eqnarray*}
Writing things in terms of quantum determinants again, we deduce
\[
\qdet{I|M}^{-1} \qdet{J|M} = \qcomm{}{J}{I} \qdet{J|M}\qdet{I|M}^{-1}\,.
\]
Finally, note that $\qdet{J|M}\qdet{I|M} = \qcomm{}{J}{I}\qdet{I|M}\qdet{J|M}$ if and only if $\qdet{J\cup M}\qdet{I\cup M} = \qcomm{}{J}{I}\qdet{I\cup M}\qdet{J\cup M}$. 
\end{proof}

We are now ready for the first advertised proof of Theorem \ref{thm:qComm}.
\begin{proof}[Proof of Theorem]
Given $J,I\subseteq[n]$ with $d=|J|\leq |I|=e$, put $s=|J\cap I|$. After Lemma \ref{thm:overlap}, we may assume $s=0$. We proceed by induction on $d$, the base case being handled in Lemma \ref{thm:disjoint}.
\smallskip

Let $j$ be the least element of $J$, i.e. $\Ell{j}{J^j}=0$, and consider $(\mathcal P_{I,J\setminus j,j})$:
\[
1=\sum_{i\in I} r_{ji}^{I\setminus i}r_{ij}^{J\setminus j}\,.
\]
In terms of quantum determinants, we have 
\[
\qdet{j|J^j}=\sum_{i\in I} \qdet{j | I^i} \qdet{i | I^i}^{-1} \qdet{i | J^j}\,.
\]
By induction, we may clear the denominator to the right and get
\begin{equation}\label{eq:JI-equals}
\qdet{j|J^j}\qdet{I}=\qcomm{}{J^j}{I^i}\sum_{i\in I} (-q)^{\Ell{i}{I^i}} \qdet{j | I^i} \qdet{i | J^j}\,.
\end{equation}
On the otherhand, we may clear the denominator on the left at the expense of $\qcomm{-}{j}{i}$:
\begin{equation}\label{eq:IJ-equals}
\qdet{I}\qdet{j|J^j}=\qcomm{-}{j}{i} \sum_{i\in I} (-q)^{\Ell{i}{I^i}} \qdet{j | I^i} \qdet{i | J^j}\,.
\end{equation}
We are nearly done. First observe the following three facts.
\[
\qcomm{}{J^j}{I^i} = \qcomm{}{J^j}{I} \qquad
\qcomm{-}{j}{i} = \qcomm{-}{j}{I} \qquad
\qcomm{}{J}{I}  =\qcomm{}{j}{I}\qcomm{}{J^j}{I}
\]
Using these observations to compare (\ref{eq:JI-equals}) and (\ref{eq:IJ-equals}) finishes the proof.
\end{proof}

\section{Preliminaries for $\mathcal G$-Proof}\label{sec:prelims-Gproof}
\subsection{Quantum flag algebra}\label{sec:TaftTowber-flags}
The algebra $\mathcal{F}\ell_q(n)$ as presented below first appeared in \cite{TafTow:1}. 

\begin{mydefi}[Quantum Flag Algebra] The quantum flag algebra $\mathcal{F}\ell_q(n)$ is the $\Bbbk$-algebra generated by symbols $\left\{ f_{I} : I\in[n]^d,\,1\leq d\leq n  \right\}$ subject to the relations indicated below.
\begin{itemize}
\item \emph{Alternating} relations $(\mathcal A_I)$: For any $I\in{[n]}^{d}$ and $\sigma\in\Sym_d$, 
\begin{equation}\label{eq:TT-Alt}
f_{I} = \left\{\begin{array}{ll}
0 & \hbox{if }I\hbox{ contains repeated indices}\\
(-q)^{-\ell(\sigma)} f_{\sigma I} & \hbox{if }\sigma I =(i_1< i_2< \cdots< i_d)
\end{array}\right.
\end{equation}

\item \emph{Young symmetry} relations $(\mathcal Y_{I,J})_{(a)}$: Fix $1\leq a\leq d\leq e\leq n-a$. For any $I\in\binom{[n]}{e+a}$ and $J\in{[n]}^{d-a}$,
\begin{equation}\label{eq:TT-YS}
0 = \sum_{\Lambda\subseteq I,|\Lambda|= a} (-q)^{-\Ell{I\setminus\Lambda}{\Lambda}} f_{I\setminus \Lambda} f_{\Lambda|J}
\end{equation}

\item \emph{Monomial straightening} relations $(\mathcal M_{J,I})$: For any $J,I\subseteq{[n]}$ with $|J|\leq |I|$,
\begin{equation}\label{eq:TT-MS}
f_{J} f_{I} = \sum_{\Lambda\subseteq I,|\Lambda|=|J|} (-q)^{\Ell{\Lambda}{I\setminus \Lambda}}f_{J|I\setminus\Lambda} f_{\Lambda}
\end{equation}
\end{itemize}
\end{mydefi}

\begin{myrema} Technically, we should have taken $I,J$ to be tuples instead of sets in (\ref{eq:TT-YS}) and (\ref{eq:TT-MS}). Identify, e.g. $I=\{i_1<i_2<\cdots < i_d\}$ with $(i_1,i_2,\ldots, i_d)$. This abuse of notation will reoccur without further ado. 
\end{myrema}

In their article, Taft and Towber construct an algebra map $\phi: \mathcal{F}\ell_q(n) \rightarrow \mathrm{M}_q(n)$ taking $f_I$ to $\qdet{I}$ and show that $\phi$ is monic, with image the subalgebra of $\mathrm{M}_q(n)$ generated by the quantum minors $\left\{\qdet{I} : I\in[n]^d,\,1\leq d\leq n \right\}$. 

We have already seen that the minors $\qdet{I}$ often $q$-commute. This relation does not appear above, and so must be a consequence of relations (\ref{eq:TT-Alt})--(\ref{eq:TT-MS}). Abbreviate the right-hand side of (\ref{eq:TT-YS}) by $Y_{I,J;(a)}$. Also, we abbreviate the difference ($lhs - rhs$) in (\ref{eq:TT-MS}) by $M_{J,I}$, and the difference ($lhs - rhs$) in (\ref{eq:qComm}) by $C_{J,I}$ (replacing $\qdet{\hbox{-}}$ by $f_{\hbox{-}}$). As (\ref{eq:qComm}),(\ref{eq:TT-YS}),(\ref{eq:TT-MS}) are all homogeneous, a likely guess is that $C_{J,I}$ is some $\Bbbk$-linear combination of a certain number of expressions $M_{K,L}$ and $Y_{M,N;(a)}$ (modulo the alternating relations). As illustrated in the example below, this simple guess works. 

\begin{myexam}[$\{1\} \surr \{2,3,4\}$] We calculate the expressions $C_{1,234}$, $M_{1,234}$, and $Y_{1234,\emptyset;(1)}$ and arrange them as rows in Table \ref{tbl:hiddenrel-example}.
Viewing the table column by column, deduce $C_{1,234} = M_{1,234} + q^2 Y_{1234,\emptyset;(1)}$.
\begin{table}[!htp]
\[
\begin{array}{|c||rrrrr|}
\hline
\rule[-1.5ex]{0pt}{2.5ex} C_{1,234} & f_{1}f_{234} &  &  &  & - q^{-1} f_{234}f_{1} \\
\hline\hline
\rule[-1ex]{0pt}{3ex} M_{1,234}& f_{1}f_{234} & - q^{2} f_{123}f_{4} & +q^{1} f_{124}f_{3} & -q^{0} f_{134}f_{2} &\\
\rule[-1.5ex]{0pt}{2.5ex} Y_{1234,\emptyset ; (1)} & & f_{123}f_{4} & -q^{-1} f_{124}f_{3} & +q^{-2} f_{134}f_{2} & -q^{-3} f_{234}f_{1}\\
\hline
\end{array}
\]
\caption{Finding the relation $f_1f_{234} - q^{-1}f_{234}f_1 = 0$.}\label{tbl:hiddenrel-example}
\end{table}
\end{myexam}

While the proof idea will be simple (``perform Gaussian elimination''), the proof itself is not. We separate out the more interesting steps below.

\subsection{POset paths}\label{sec:poset-paths}
Given a set $X$, the elements of the power set $\mathcal P X$ have a partial ordering: for $A,B\in\mathcal P X$, we say $A<B$ if $A\subsetneq B$. We are interested in the case $X\subseteq[n]$ and we think of this POset as an edge-weighted, directed graph as follows.

\begin{mydefi} Given $I,J\subseteq[n]$ such that $J\surr I$, the graph $\Gamma(J;I)$ has vertex set $\mathcal{V}=\mathcal PJ$ and edge set $\mathcal{E}=\{(A,B) \mid A,B\in \mathcal{V},\,A\subsetneq B\}$. Each edge $(A,B)$ of $\Gamma$ has a weight $\alpha_A^B$ given by the function $\alpha:\mathcal{E} \rightarrow \Bbbk$,
\begin{equation}\label{eq:alpha-def}
\forall (A,B)\in\mathcal{E}: \quad \alpha_{A}^B = (-q)^{-\Ell{J\setminus B}{{B} \setminus A} - \Ell{B\setminus A}{A} +\big(2|J\setminus B|-|I|\big)|(B\setminus A)\cap J'|}
\end{equation}
for $J'$ is as in Definition \ref{def:surrounds}.
\end{mydefi}

\begin{myexam} If $|J|=m$, then $\Gamma(J)$ has $2^m$ vertices and $\sum_{k=1}^m\binom{m}{k}(2^m-1)$ edges. In Figure \ref{fig:gamma-156}, we give an illustration of $\Gamma(\{1,5,6\})$, omitting two edges and many edge weights for legibility.

\begin{figure}[!htbp]
$$
\xymatrix@!0{
& 56 \ar@{->}[rr] 
& & 156 
\\
6 \ar@{->}[ur]\ar@{->}[rr]\ar@{->}[urrr]
& & 16 \ar@{->}[ur]
\\
& 5 \ar@{->}[rr]\ar@{->}[uu]
& & 15\ar@{->}[uu]_{\alpha_{15}^{156}}
\\
\emptyset \ar@{->}[uu]^{\alpha_{\emptyset}^{6}} \ar@{->}[urrr]\ar@{->}[rr] \ar@{->}[ur] \ar@{->}[uuur]
& & 1 \ar@{->}[ur]_{\alpha_{1}^{15}} \ar@{->}[uuur]\ar@{->}[uu]
}
$$
\caption{The graph $\Gamma(\{1,5,6\})$ (partially rendered).} \label{fig:gamma-156}
\end{figure}
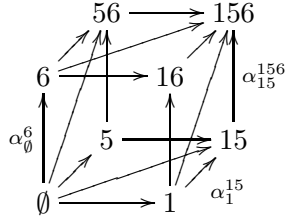
\end{myexam}

For the remainder of the subsection, we assume $J\cap I=\emptyset$. Write $J=J'\dot{\cup}J'' = \{j_1<\ldots<j_{r'}\}\cup\{j_{r'+1}<\cdots<j_{r'+r''}\}$; also, put $|J|=r'+r''=r$, $|I|=s$, and $s-r=t$. 

In the graph $\Gamma(J;I)$, we consider \emph{paths} and \emph{path weights} defined as follows:
\[
\mathfrak{P}_0 = \Big\{(A_1,A_2,\ldots,A_p) \mid A_i\subseteq J \hbox{ s.t. }\emptyset\subsetneq A_1\subsetneq A_2\subsetneq\cdots\subsetneq A_p\subsetneq J \Big\}
\]
and $\mathfrak{P}=\mathfrak{P}_0 \cup \hat{0}\cup \hat{1}$, where $\hat{0}=(\emptyset)$, and 
\[
\hat{1} = (\{j_{r'+1}\},\{j_{r'+1},j_{r'+2}\},\ldots,J'',\{j_{r'},\ldots, j_{r}\},\ldots,\{j_2,\ldots, j_r\},J).
\]
The \emph{weight} $\alpha(\pi)$ of a path $\pi=(A_1,\ldots,A_p)\in\mathfrak{P}_0$ is the product of edge weights of the augmented path $(\emptyset,\pi,J)$:
\[
\alpha_{\emptyset}^{A_1}\cdot\alpha_{A_1}^{A_2}\cdots \alpha_{A_{p-1}}^{A_p}\cdot\alpha_{A_p}^{J}.
\]
We extend the definition of $\alpha$ to all of $\mathfrak{P}$ as follows. Notice that if $B=A$ in (\ref{eq:alpha-def}), we get $\alpha_{A}^{A}=1$. With this broader definition of the weight function $\alpha$, we may define $\alpha(\pi) = \alpha(\emptyset,\pi,J)$ for  $\pi=\hat{0},\hat{1}$ as well. Writing $\hat{1}= (A_1,\ldots,A_{r=|J|})$, the path $(A_1,\ldots,A_{r-1})\in\mathfrak{P}_0$ will also be important. We label this special path $\pi^{\hat{1}}$.

\begin{mydefi} Given a subset $K\subseteq J$, define $\mM(K)$ as follows. If $ K\cap J' \neq \emptyset$, put $\mM(K) = \min(K\cap J')$. Otherwise, put $\mM(K) = \max(K\cap J'')$. 
\end{mydefi}

For any path $\pi=(A_1,\ldots,A_p)$, put $A_0=\emptyset$ and $A_{p+1}=J$. Notice that $\hat{1}$ has the property that $A_{k}\setminus A_{k-1} \neq \mM(A_{k+1}\setminus A_{k-1})$ for all $1\leq k< r$, but $A_r=\mM(A_{r+1}\setminus A_{r-1})$.

\begin{mydefi} Fix a length $1\leq p\leq r-1$. A path $(A_1,\ldots,A_p)\in\mathfrak{P}_0$ shall be called \emph{regular} (or \emph{regular at position $i_0$}), if $(\exists i_0) (1\leq i_0\leq p)$ satisfying: (a) $|A_i |= i\,(\forall\, 1\leq i\leq i_0) $; (b) $A_{i_0} \setminus A_{i_0-1} = \mM(A_{i_0+1} \setminus A_{i_0-1} )$ (again, taking $A_0=\emptyset$ and $A_{p+1}=J$ if necessary). A sequence is called \emph{irregular} if it is nowhere regular. Extend the notion of regularity to $\mathfrak{P}$ by calling $\hat{0}$ irregular and $\hat{1}$ regular.
\end{mydefi}

\begin{myrema} The set $\mathfrak{P}$ is the disjoint union of its regular and irregular paths. We point out this tautology only to emphasize its importance in the coming proposition. Write $\mathfrak{P}'$ for the irregular paths, and $\mathfrak{P}''$ for the regular paths. \end{myrema}

\begin{myprop}\label{thm:paths-bijection} The subsets $\mathfrak{P}'$ and $\mathfrak{P}''$ of $\mathfrak{P}$ are equinumerous.
\end{myprop}

We will build a bijective map $\wp$ between the two sets. Given an irregular path $\pi=(A_1,\ldots,A_p)\in\mathfrak{P}_0$, we insert a new set $B$ so that $\wp(\pi)$ is regular at $B$:
\begin{enumerate}
\item Find the unique $i_0$ satisfying: $(|A_i|=i\quad \forall i\leq i_0)\wedge(|A_{i_0+1}|>i_0+1)$.
\item Compute $b=\mM(A_{i_0+1}\setminus A_{i_0})$ 
\item Put $B=A_{i_0} \cup \{b\}$.
\item Define $\wp(\pi) := (A_1,\ldots,A_{i_0},B,A_{i_0+1},\ldots,A_p)$.
\end{enumerate}
For the remaining irregular path $\hat{0}$, we put $\wp(\hat{0}) = (\{j_1\})$, which agrees with the general definition of $\wp$ if we think of $\hat{0}$ as the empty path $()$ instead of the path consisting of the empty set.

\begin{myexam} Table \ref{tbl:varphi-example} illustrates the action of $\wp$ on  $\mathfrak{P}$ when $J=\{1,5,6\}$.
\begin{table}[!htp]
\[
\begin{array}{|@{\;\;}c@{\;\;}||@{\;\;}c@{\;\;}|@{\;\;}c@{\;\;}|@{\;\;}c@{\;\;}|@{\;\;}c@{\;\;}|@{\;\;}c@{\;\;}|@{\;\;}c@{\;\;}|@{\;\;}c@{\;\;}|} 
\hline
\rule[-1.4ex]{0pt}{4ex}\pi & \hat 0 & (5) & (6) & (15) & (16) & (56) & (5,56) \\
\hline
\rule[-1.4ex]{0pt}{4ex}\wp(\pi) & (1)  & (5,15)  & (6,16)  & (1,15)  & (1,16)  & (6,56)  & \hat 1\\
\hline
\end{array}
\]
\caption{The pairing of $\mathfrak{P}'$ and $\mathfrak{P}''$ via $\wp$.}\label{tbl:varphi-example}
\end{table}
\end{myexam}

\begin{proof}[Proof of Proposition] We reach a proof in three steps.
\par\noindent\textit{Claim 1:} $\wp(\mathfrak P')\subseteq \mathfrak P''$.
\par Take a path $\pi\in\mathfrak P'$ (i.e. a path with no regular points). The effect of $\wp$ is to insert a regular point at position $i_0+1$ (the spot where $B$ sits), so the claim is proven if we can show $\wp(\pi)\in\mathfrak{P}$. 

As $\wp(\hat{0})$ clearly belongs to $\mathfrak{P}$, we may focus on those $\pi\in\mathfrak{P}_0$. Also, it is plain to see that $\pi^{\hat{1}}$ is irregular, and $\wp(\pi^{\hat{1}})=\hat{1}$. If $\wp$ is to be a bijection, we are left with the task of showing that $\wp(\mathfrak{P}'\cap\mathfrak{P}_0\setminus \pi^{\hat{1}})\subseteq\mathfrak{P}_0$

When $|A_p|<r-1$, any $B$ that is inserted will result in another path in $\mathfrak{P}_0$ (because $|B|$ must be less than $r$). When $|A_p|=r-1$, there is some concern that we will have to insert a $B$ at the end of the path, resulting in $J$ being the new terminal vertex---disallowed in $\mathfrak{P}_0$. This cannot happen:
\par\textit{Case $p<r-1$:} At some point $1\leq i_0<p$, there is a jump in set-size greater than one when moving from $A_{i_0}$ to $A_{i_0+1}$. Hence, the $B$ to be inserted will not come at the end, but rather immediately after $A_{i_0}$ to $A_{i_0+1}$

\par\textit{Case $p=r-1$:} The only path $(A_1,A_2,\ldots,A_{r-1})\in\mathfrak{P}_0$ which is nowhere regular is the path $\pi^{\hat{1}}$.
\medskip

\par\noindent\textit{Claim 2:} $\wp$ is 1-1.
\par Suppose $\wp(A_1,\ldots,A_p) = \wp(A'_1,\ldots, A'_{p'})$, and suppose we insert $B$ and $B'$ respectively. By the nature of $\wp$, we have $p=p'$ and $i_0 \neq i'_0$. Take $i_0<i'_0$.
Also notice that $(A'_1,\ldots,A'_{p'}) = (A_1,\ldots,A_{i_0},B, A_{i_0+1},\ldots,A'_{i'_0}, \ldots A'_{p'})$ In particular, $B$ is a regular point of $(A'_1,\ldots,A'_{p'})$, and consequently, \\$(A'_1,\ldots,A'_{p'})\not\in\mathfrak{P}'$. 
\medskip

\par\noindent\textit{Claim 3:} $\wp$ is onto.
\par Consider a path $\pi=(A_1,\ldots,A_p)\in\mathfrak{P}''$. If $p=1$, then it is plain to see that the only irregular path is $\pi= (\{j_1\})$, which is the image of $(\emptyset)$ under $\wp$. So we consider $\pi\in\mathfrak{P}''$ with $p>1$. Note that $|A_1| = 1$, for otherwise $\pi$ cannot have any regular points. Now, locate the first $1\leq i_0\leq p$ with (a) $|A_{i_0}| = i_0$; and (b) $A_{i_0}\setminus A_{i_0-1} = \mM(A_{i_0+1}\setminus A_{i_0-1}$. The path $\pi'=(A_1,\ldots,A_{i_0-1},A_{i_0+1},\ldots,A_k)$ is in $\mathfrak{P}'$ and moreover, $\wp(\pi')=\pi$. 
\end{proof}

Certainly one could cook up other bijections between the regular and irregular paths in $\mathfrak{P}$. The map we have used has an additional nice property. 

\begin{myprop} The bijection $\wp$ from the proof of Proposition \ref{thm:paths-bijection} is path-weight preserving.\end{myprop}

The result rests on 
\begin{mylemm}\label{paths:weight-preserving} Let $\emptyset\subseteq A\subseteq B\subseteq C\subseteq J$. Writing $\hat{B}=B\setminus A$ and $\hat{C}=C\setminus B$, we have
\begin{equation}\label{eq:alpha-two-to-one}
\alpha_A^B \alpha_B^C = \left[(-q)^{2\Ell{B'\cap J' }{C'} - 2\Ell{C' }{B'\cap J''}} \right]\alpha_A^C\,.
\end{equation}
\end{mylemm}

\begin{proof}
From the definition of $\alpha_B^C$, we have
\begin{eqnarray*}
\alpha_A^B &=& (-q)^{-\Ell{J\setminus B }{\hat{B}} -\Ell{\hat{B}}{A} +\big(2|J\setminus B|-|I|\big)|\hat{B}\cap J'|} \\
\alpha_B^C &=& (-q)^{-\Ell{J\setminus C }{\hat{C}} -\Ell{\hat{C}}{B} +\big(2|J\setminus C|-|I|\big)|\hat{C}\cap J'|} \\
\alpha_A^C &=& (-q)^{-\Ell{J\setminus C }{\hat{B}\cup\hat{C}} -\Ell{\hat{B}\cup\hat{C}}{A} +\big(2|J\setminus C|-|I|\big)|(\hat{B}\cup\hat{C})\cap J'|} 
\end{eqnarray*}
Let us compare the exponents of $\alpha_A^C$ and $\alpha_A^B\alpha_B^C $:
\begin{eqnarray}
\nonumber \exp(\alpha_A^C) &=& -\Ell{J\setminus C }{\hat{B}}-\Ell{J\setminus C }{\hat{C}} -\Ell{\hat{C}}{A} - \Ell{\hat{B}}{A} +\\
\label{eq:one-alpha} && \big(2|J\setminus A|-2|\hat{C}|-2|\hat{B}|-|I|\big)\big(|\hat{B}\cap J'| + |\hat{C}\cap J'|\big)\,,
\end{eqnarray}
while 

\begin{eqnarray}
\nonumber \exp(\alpha_A^B\alpha_B^C ) &=& -\Ell{J\setminus B }{\hat{B}} - \Ell{J\setminus C }{\hat{C}} -\Ell{\hat{B}}{A} -\Ell{\hat{C}}{B} +\\
\nonumber && \big(2|J\setminus B|-|I|\big)|\hat{B}\cap J'| +\big(2|J\setminus C|-|I|\big)|\hat{C}\cap J'| \\
\nonumber &=& -\Big\{\Ell{J\setminus C}{\hat{B}}+\Ell{\hat{C}}{\hat{B}}\Big\} - \Ell{J\setminus C}{\hat{C}} -\Ell{\hat{B}}{A} - \\
\nonumber && \Big\{\Ell{\hat{C}}{A}+\Ell{\hat{C}}{\hat{B}}\Big\} + \Big\{2|J\setminus A|-2|\hat{B}|-|I|\Big\}|\hat{B}\cap J'| + \\
\nonumber &&  \Big\{2|J\setminus A| -2|\hat{B}|-2|\hat{C}|-|I|\Big\}|\hat{C}\cap J'| \\
\label{eq:two-alphas} &=& 2|\hat{C}||\hat{B}\cap J'|- 2\Ell{\hat{C}}{\hat{B}} +  \Big\{\exp(\alpha_A^C)\Big\}\,.
\end{eqnarray}
Notice that $2|\hat{C}||\hat{B}\cap J'|=2\Ell{\hat{C}}{\hat{B}\cap J'} + 2\Ell{\hat{B}\cap J'}{\hat{C}}$, and that $-2\Ell{\hat{C}}{\hat{B}} = -2\Ell{\hat{C}}{\hat{B}\cap J'} -2\Ell{\hat{C}}{\hat{B}\cap J''}$. The discrepancy between (\ref{eq:one-alpha}) and (\ref{eq:two-alphas}) becomes
$2\Ell{\hat{B}\cap J'}{\hat{C}}  -2\Ell{\hat{C}}{\hat{B}\cap J''}$, as desired.
\end{proof}

Now the proposition follows by comparing $\alpha(A_{i_0},A_{i_0+1})$ and $\alpha(A_{i_0},B,A_{i_0+1})$.

\begin{proof}[Proof of Proposition] Suppose that $\pi=(\ldots,A,C,\ldots)$, and that $\wp(\pi)$ inserts $B$ immediately after $A$. Then $B=A\cup\mM(C\setminus A)$. Writing $b=\mM(C\setminus A)$, (\ref{eq:alpha-two-to-one}) implies
$$
\alpha(\wp(\pi)) = \left[(-q)^{2\Ell{b\cap J'}{\hat{C}}  -2\Ell{\hat{C}}{b\cap J''}} \right]\cdot\alpha(\pi)\,.
$$
Now, if $b\cap J'\neq\emptyset$, then $b$ is the smallest element in $C\setminus A$, and in particular, $\Ell{b}{\hat{C}}=0$. In this same case, $b\cap J''=\emptyset$, so $\Ell{\hat{C}}{b\cap J''} =0$ too. An analogous argument works for the case $b\cap J'=\emptyset$.
\end{proof}

One more interesting fact about $\Gamma(J;I)$ and $\mathfrak{P}$ is worth mentioning. When calculating $\alpha(\pi^{\hat{1}})$ using (\ref{eq:alpha-two-to-one}), the twos introduced in the exponents there all disappear.

\begin{myprop} Given, $J, J',J''$, and $\pi^{\hat{1}}$ as above, we have
\begin{equation}\label{eq:alpha-of-special-path}
\alpha(\pi^{\hat{1}}) = (-q)^{|J'|(|J'|-1)-|J''|(|J''|-1)}\times\alpha_{\emptyset}^J\,.
\end{equation}
\end{myprop}

\begin{proof} Applying (\ref{eq:alpha-two-to-one}) repeatedly to the expression $\alpha(\pi^{\hat{1}})$ we see that
\begin{eqnarray*}
\alpha(\pi^{\hat{1}}) &=& \left[(-q)^{2\Ell{j_{r'+1}\cap J' }{j_{r'+2}} - 2\Ell{j_{r'+2}}{ j_{r'+1}\cap J''}}\right] \times \\
&& \alpha_{\emptyset}^{j_{r'+1}j_{r'+2}} \alpha_{j_{r'+1}j_{r'+2}}^{j_{r'+1}j_{r'+2}j_{r'+3}} \cdots \alpha_{j_{2}\cdots j_{r}}^J\\
 &=& (-q)^{-2(1)}\left[(-q)^{2\Ell{j_{r'+2}\cap J'}{ j_{r'+3}} - 2\Ell{j_{r'+3} }{j_{r'+2}\cap J''}}\right] \times \\
 && \alpha_{\emptyset}^{j_{r'+1}j_{r'+2}j_{r'+3}} \cdots \alpha_{j_{2}\cdots j_{r}}^J\\
 &=& (-q)^{-2(1)-2(2)}\left[(-q)^{2\Ell{j_{r'+3}\cap J'}{ j_{r'+4}} - 2\Ell{j_{r'+4}}{ j_{r'+3}\cap J''}}\right] \times \\
&& \alpha_{\emptyset}^{j_{r'+1}j_{r'+2}j_{r'+3}j_{r'+4}} \cdots \alpha_{j_{2}\cdots j_{r}}^J\\
&\vdots& \\
 &=& (-q)^{-2(1)-\cdots -2(|J''|-1)}\left[(-q)^{2\Ell{j_{r}\cap J' }{ j_{r'}} - 2\Ell{j_{r'}}{ j_{r}\cap J''}}\right] \times \\
 && \alpha_{\emptyset}^{j_{r'}\cdots j_{r}} \cdots \alpha_{j_{2}\cdots j_{r}}^J\\
&=& (-q)^{-2\frac{(|J''|-1)|J''|}{2}}(-q)^{0-0}\left[(-q)^{2\Ell{j_{r'}\cap J' }{ j_{r'-1}}- 2\Ell{j_{r'-1}}{ j_{r'}\cap J''}}\right] \times \\
&& \alpha_{\emptyset}^{j_{r'-1}\cdots j_{r}} \cdots \alpha_{j_{2}\cdots j_{r}}^J\\
&=& (-q)^{2(1)} (-q)^{-|J''|(|J''|-1)}\left[(-q)^{2\Ell{j_{r'-1}\cap J'}{ j_{r'-2}} - 2\Ell{j_{r'-2}}{ j_{r'-1}\cap J''}}\right] \times \\
&& \alpha_{\emptyset}^{j_{r'-2}\cdots j_{r}} \cdots \alpha_{j_{2}\cdots j_{r}}^J\\
&\vdots& \\
&=& (-q)^{2(1)+\cdots +2(|J'|-1)}(-q)^{-|J''|(|J''|-1)}\times\alpha_{\emptyset}^J \\
&=& (-q)^{|J'|(|J'|-1)-|J''|(|J''|-1)}\times\alpha_{\emptyset}^J \,. \,\, \qedhere
\end{eqnarray*}
\end{proof}

\def\MessyNotation{Only minor changes to this proof are needed to prove the theorem in the general setting (e.g. replacing every instance of $J$ below with $J_0:=J\setminus I$). In the interest of avoiding even more notation, we leave this work to the reader.}
\section{$\mathcal G$-Proof of Theorem}\label{sec:Gproof}
We keep the notations $J',J'', r',r'',r,s,t$ from Section \ref{sec:poset-paths}, and as we did there, we only consider the case $J\cap I =\emptyset$.\footnote{\MessyNotation} Before we dive in, we define a new quantity $CM_{J,I}(\theta)$. 
\begin{eqnarray*}
C_{J,I} - M_{J,I} &=& - q^{|J''|-|J'|} f_I f_J + \left(\sum_{\Lambda\subseteq I,\,|\Lambda|=r} (-q)^{\Ell{\Lambda}{I^\Lambda}} f_{J|I\setminus \Lambda} f_{\Lambda} \right) \\
&=& \sum_{\Lambda\subseteq I} (-q)^{|J'|t}(-q)^{-\Ell{I^\Lambda}{\Lambda}} f_{J \cup(I\setminus \Lambda)}f_{\Lambda}  - q^{|J''|-|J'|} f_I f_J  \\
&=& \sum_{\Lambda\subseteq I} (-q)^{|J'|t+|J''||J|}(-q)^{-\Ell{(J\cup I)^\Lambda}{\Lambda}} f_{(J \cup I)\setminus \Lambda} f_{\Lambda}  - q^{|J''|-|J'|} f_I f_J  
\end{eqnarray*}
Here, we have replaced $\Ell{\Lambda}{I^\Lambda}$ with $|I^\Lambda||\Lambda|-\Ell{I^\Lambda}{\Lambda}$ and $\Ell{J}{I^\Lambda}$ with $|J||I^\Lambda|-\Ell{I^\Lambda}{J}$. 
\[
CM_{J,I}(\theta) :=  (-q)^{|J'|t+|J''||J|}\left(\sum_{\Lambda\subseteq I} (-q)^{-\Ell{(J\cup I)^\Lambda}{\Lambda}} f_{(J \cup I)\setminus \Lambda} f_{\Lambda}  - \theta f_I f_J\right) \,.
\]
We prove the theorem in steps:

\begin{myprop}\label{thm:Comm-via-YS-qComm-I} Suppose $I,J\subseteq[n]$ are such that $J\surr I$. With $CM_{J,I}(\theta)$ and $Y_{L,K;(a)}$ as defined above, 
\[
\sum_{\emptyset\subseteq K\subsetneq J} \eta_K \cdot Y_{(I\cup J)\setminus K,K;(r-|K|)} = CM_{J,I}(\theta) 
\]
for some constants $\{\eta_K \in Z[q,q^{-1}]: \emptyset\subseteq K\subsetneq J\}$ and $\theta\in Z[q,q^{-1}]$.
\end{myprop}

\begin{myprop}\label{thm:Comm-via-YS-qComm-II} In the notation above, $\theta = (-q)^{-|J'|t-|J''||J|}q^{|J''|-|J'|}$.
\end{myprop}

The alternating property of the symbols $f_K$ and the product in $\mathcal{F}\ell_q(n)$ play no role in our proof, so we begin by eliminating these distractions. Let $V$ be the vector space over $\Bbbk$ with basis $\{e_{(A,B)} : A\cup B=I\cup J, A\cap B=\emptyset,\hbox{ and }|B|=r\}$. There is a $\Bbbk$-linear map $\mu:V\rightarrow \mathcal{F}\ell_q(n)$, sending $e_{A,B}$ to $f_A f_B$. The vectors
\[
v^{\theta} := \sum_{\Lambda\subseteq I} (-q)^{-\Ell{(I\cup J)^\Lambda}{\Lambda}} e_{(I \cup J)\setminus\Lambda,\Lambda}  - \theta e_{I,J}
\]
and (for each $\emptyset\subseteq K\subsetneq J$) 
\[
v^K := \sum_{\Lambda\subseteq(I\cup J),|\Lambda|=r-|K|} (-q)^{-\Ell{(I\cup J^K)^\Lambda}{\Lambda}}(-q)^{-\Ell{\Lambda}{K}} e_{(I\cup J)\setminus (K\cup \Lambda),K\cup \Lambda}
\]
have familiar images. Check that $\mu((-q)^{|J'|t+|J''||J|}\cdot v^{\theta})= CM_{J,I}(\theta)$ and $\mu(v^K) = Y_{(I\cup J)\setminus K,K;(r-|K|)}$. 

Proposition \ref{thm:Comm-via-YS-qComm-I} will be proven if we can show that $v^{\theta}$ is a linear combination of the $v^K$ for some $\theta$. This is not immediate as the span of the vectors $v^K$ has dimension (at most, \emph{a priori}) $2^{r}-1$, while $V$ is $\binom{r+s}{r}$ dimensional. 

\begin{mydefi} For each $K\in\mathcal PJ$, let $V_{(K)} = \mathrm{span}_{\Bbbk}\big\{ e_{A,B} : B\cap J = K \big\}$. Clearly, $V$ is graded by the POset $\mathcal PJ$, i.e., $V=\bigoplus_{K\in\mathcal PJ} V_{(K)}$. For each $K\in\mathcal PJ$, define the distinguished element $e^{K}$ by
\[
e^K = \sum_{\Lambda\subseteq I, |\Lambda|=r-|K|} (-q)^{-\Ell{(I\cup J)^{(K\cup \Lambda)}}{\Lambda}}(-q)^{-\Ell{\Lambda}{K}}e_{(I\cup J)\setminus(\Lambda\cup K),\Lambda\cup K}.
\]
For any $v\in V$, write $(v)_{(K)}$ for the component of $v$ in $V_{(K)}$, that is, $v=\sum_{K} (v)_{(K)}$.
\end{mydefi}

Notice that $e^J = e_{I,J}$, and that
\[
e^{\emptyset} = \sum_{\Lambda\subseteq I,|\Lambda|=r} (-q)^{-\Ell{(I\cup J)^\Lambda}{\Lambda}}e_{(I\cup J)\setminus\Lambda,\Lambda}
\]
In other words, ${v}^{\theta} = e^{\emptyset}-\theta e^J$. Good fortune provides that the $v^{K'}$ may also be expressed in terms of the $e^K$.

\begin{mylemm} For each $K'\in\mathcal PJ\setminus J$, there are constants $\alpha_{K'}^{K}\in \Bbbk$ satisfying
\[
v^{K'} = \sum_{K\in\mathcal PJ} \alpha_{K'}^{K} e^K\,.
\]
\end{mylemm}

\begin{myrema} As the proof will show, these $\alpha_{K'}^K$ are precisely the edge-weights of $\Gamma(J;I)$ from Section \ref{sec:poset-paths}, in particular $\alpha_{K}^K = 1$. It will also show that $\alpha_{K'}^{K}=0$ if $K'\nless K$ in the POset $\mathcal PJ$, a critical ingredient in the approaching Gaussian elimination argument. 
\end{myrema}

\begin{proof}[Proof of Lemma] Fixing a subset $K'$, if $K\supsetneq K'$, we write $\hat{K} = K\setminus K'$. Similarly, let $\hat{\Lambda} = \Lambda\setminus J$. Studying $v^{K'}$, we see that
\begin{eqnarray*}
v^{K'} &=& \sum_{\newatop{\Lambda\subseteq(I\cup J)\setminus K'}{|\Lambda|=r-|K'|}} (-q)^{-\Ell{(I\cup J^{K'})^{\Lambda}}{\Lambda}}(-q)^{-\Ell{\Lambda}{K'}} e_{(I\cup J)\setminus (\Lambda\cup K'),\Lambda\cup K'} \\
&=& \sum_{K\in\mathcal PJ} (v^{K'})_{(K)} \\
&=& \sum_{K\in\mathcal PJ} \sum_{\newatop{\Lambda\subseteq(I\cup J)\setminus K'}{ \Lambda\cap J = \hat{K}}} (-q)^{-\Ell{(I\cup J)^{(\hat{\Lambda}\cup K)}}{\hat{\Lambda}\cup\hat{K}}} \times \\
&& \rule[-3.5ex]{0pt}{7ex} (-q)^{-\Ell{\hat{\Lambda}\cup\hat{K}}{K'}} e_{(I\cup J)\setminus (\hat{\Lambda}\cup K),\hat{\Lambda}\cup K} 
\end{eqnarray*}
\begin{eqnarray*}
&=& \sum_{K\in\mathcal PJ} (-q)^{-\Ell{(I^{\hat{\Lambda}}) \cup (J^K)}{\hat{K}}}(-q)^{-\Ell{\hat{K}}{K'}} \times \\
&& \left(\sum_{\newatop{\hat{\Lambda}\subseteq I}{ |\hat{\Lambda}|=r-|K|}} (-q)^{-\Ell{(I\cup J)^{(\hat{\Lambda}\cup K)}}{\hat{\Lambda}}}(-q)^{-\Ell{\hat{\Lambda}}{K'}} e_{(I\cup J)\setminus (\hat{\Lambda}\cup K),\hat{\Lambda}\cup K} \right).
\end{eqnarray*}
Why can we perform this last step? Because $J\surr I$, the expression $\Ell{I^{\hat{\Lambda}}}{\hat{K}}$ does not actually depend on $\hat{\Lambda}$, only on $|\hat{\Lambda}|$. Indeed, it equals $|I\setminus{\hat{\Lambda}}|\cdot|\hat{K}\cap J'|$. Multiplying and dividing by $(-q)^{-\Ell{\hat{\Lambda}}{\hat{K}}}$, we rewrite this last expression as
\begin{eqnarray*}
v^{K'} &=& \sum_{K} (-q)^{-\Ell{(I^{\hat{\Lambda}}) \cup (J^K)}{\hat{K}}}(-q)^{-\Ell{\hat{K}}{K'}+\Ell{\hat{\Lambda}}{\hat{K}}} \times \\
&& \left(\sum_{\hat{\Lambda}\subseteq I, |\hat{\Lambda}|=r-|K|} (-q)^{-\Ell{(I\cup J)^{(\hat{\Lambda}\cup K)}}{\hat{\Lambda}}}(-q)^{-\Ell{\hat{\Lambda}}{K}} e_{(I\cup J)\setminus (\hat{\Lambda}\cup K),\hat{\Lambda}\cup K} \right) \\
&=& \sum_{K'\leq K}(-q)^{\big(2|J\setminus K|-|I|\big)|\hat{K}\cap J'| - \Ell{J^K}{\hat{K}} -\Ell{\hat{K}}{K'}} \times \left(e^K\right)\\
&=& \sum_{K'\leq K} \alpha_{K'}^K e^K\,.\qedhere
\end{eqnarray*}
\end{proof}

\begin{mycoro} For any $v^{K'},v^K$ with $K'<K$ in the POset $\mathcal PJ$, and for the same constants $\alpha_{K'}^K$ as defined above, we have
$$
(v^{K'}-\alpha_{K'}^{K} v^K)_{(K)} = 0.
$$
\end{mycoro}

\begin{proof}[Proof of Proposition \ref{thm:Comm-via-YS-qComm-I}] We use the corollary to perform a certain Gaussian elimination on the ``matrix'' of the vectors $v^K$. Table \ref{tbl:vK-eK-coords} displays this matrix for the POset $\mathcal P(\{1,5,6\})$ should make our intentions clear. 
\medskip
\begin{table}[!htp]
$$
\begin{array}{c||c|ccc|ccc|c}
& e^{\emptyset} & e^1 & e^5 & e^6 & e^{15} & e^{16} & e^{56} & e^{156} \\
\hline
\hline
\rule[-1.5ex]{0pt}{4ex} 
v^{15} & & & & & 1 & & & \alpha_{15}^{156} \\
\rule[-1.5ex]{0pt}{4ex} 
v^{16} & & & & & & 1 & & \alpha_{16}^{156} \\
\rule[-1.5ex]{0pt}{4ex} 
v^{56} & & & & & & & 1 & \alpha_{56}^{156} \\
\hline
\rule[-1.5ex]{0pt}{4ex} 
v^{1} &  & 1 & & & \alpha_{1}^{15} & \alpha_{1}^{16} & & \alpha_{1}^{156} \\
\rule[-1.5ex]{0pt}{4ex} 
v^{5} &  & & 1 & & \alpha_{5}^{15} & & \alpha_{5}^{56} & \alpha_{5}^{156} \\
\rule[-1.5ex]{0pt}{4ex} 
v^{6} &  & & & 1 & & \alpha_{6}^{16} & \alpha_{6}^{56} & \alpha_{6}^{156} \\
\hline
\rule[-1.5ex]{0pt}{4ex} 
v^{\emptyset} & 1 & \alpha_{\emptyset}^{1} & \alpha_{\emptyset}^{5} & \alpha_{\emptyset}^{6} & \alpha_{\emptyset}^{15} & \alpha_{\emptyset}^{16} & \alpha_{\emptyset}^{56} & \alpha_{\emptyset}^{156} 
\end{array} 
$$ 
\caption{Writing the vectors $v^{K'}$ in terms of the $e^K$.}\label{tbl:vK-eK-coords}
\end{table}
\medskip

Performing Gaussian elimination between the rows in the first two layers of the matrix, we see that the new rows in the second layer---who began their life with $|J|+1$ nonzero entries---now have exactly two nonzero entries.
\begin{eqnarray*}
(v^{J\setminus\{k,l\}})' &=& v^{J\setminus\{k,l\}} - \alpha_{J\setminus \{k,l\}}^{J\setminus k} v^{J\setminus k} - \alpha_{J\setminus \{k,l\}}^{J\setminus l} v^{J\setminus l} \\
&=& e^{J\setminus \{k,l\}} + \left(\alpha_{J\setminus \{k,l\}}^J - \alpha_{J\setminus \{k,l\}}^{J\setminus k}\alpha_{J\setminus k}^{J} -  \alpha_{J\setminus \{k,l\}}^{J\setminus l}\alpha_{J\setminus l}^{J} \right)e^J\,,
\end{eqnarray*}
e.g., $v^1$ from Table \ref{tbl:vK-eK-coords} becomes $(v^1)' = e^{1} + \left(\alpha_{1}^{156} - \alpha_{1}^{15}\alpha_{15}^{156} -  \alpha_{1}^{16}\alpha_{16}^{156} \right)e^{156}$. 
Marching down the layers of this matrix one-by-one, we see that the new final row is given by
$(v^{\emptyset})' = e^{\emptyset} + \theta e^{J}={v}^{\theta}$ for some $\theta$.
\end{proof}

\begin{proof}[Proof of Proposition \ref{thm:Comm-via-YS-qComm-II}] Careful bookkeeping shows that
\begin{eqnarray}
\nonumber \theta &=& \alpha_{\emptyset}^{J} -\left(\sum_{\emptyset\subsetneq K\subsetneq J} \alpha_{\emptyset}^{K}\alpha_{K}^{J}\right) + \left(\sum_{\emptyset\subsetneq K_1\subsetneq K_2\subsetneq J} \alpha_{\emptyset}^{K_1} \alpha_{K_1}^{K_2} \alpha_{K_2}^{J}\right) -\cdots \\
&& \cdots + (-1)^{|J|-1}\left(\sum_{\emptyset\subsetneq K_1\subsetneq\cdots\subsetneq K_{|J|-1}\subsetneq J} \alpha_{\emptyset}^{K_1} \alpha_{K_1}^{K_2} \cdots \alpha_{K_{r-1}}^{J}\right)\,.
\end{eqnarray}
In other words, $\theta$ is a signed sum of path weights $\alpha(\pi)$, $\pi$ running over all paths in $\mathfrak{P}$ save for $\hat{1}$. As the sign attached to $\pi$ is the same as the length of $\pi$, and as the bijection $\wp$ from Section \ref{sec:poset-paths} increases length by one but preserves path weight, we immediately conclude
\begin{eqnarray*}
\theta &=& (-1)^{|J|-1}\alpha(\pi^{\hat{1}}) \\
&=& (-1)^{|J|-1}(-q)^{|J'|(|J'|-1)-|J''|(|J''|-1)}\times\alpha_{\emptyset}^J \\
&=& (-1)^{|J|-1}(-q)^{|J''|-|J'|}(-q)^{|J'||J'| -|J''||J''| -|I||J'|} \\
&=& q^{|J''|-|J'|}(-q)^{|J'||J'| -|J''||J''| -|J''||J'| -|I||J'| - (|J'|+|J''|+t)|J'|+|J''||J'| } \\
&=& q^{|J''|-|J'|}(-q)^{-|J'|t-|J''||J|}\,. \qedhere
\end{eqnarray*}
\end{proof}

With Proposition \ref{thm:Comm-via-YS-qComm-II} proven, Theorem \ref{thm:qComm} is finally demonstrated (modulo the Taft-Towber isomorphism $\phi$). Moreover, we achieve the second goal stated in the introduction. A brief discussion of the first goal follows.

\section{On Quantum- and Quasi- Flag Varieties}
The algebra $\mathcal{F}\ell_q(n)$ is a quantum deformation of the classic multihomogeneous coordinate ring of the full flag variety over $\mathrm{GL}_n$. The deformation was constructed in a somewhat ad-hoc manner, and we would like to know whether a theory of \emph{noncommutative flag varieties} using quasideterminants could help explain the choices for the relations in $\mathcal{F}\ell_q(n)$. In \cite{Lau:1}, it is shown that any relation $(\mathcal{Y}_{I,J})_{(a)}$ has a quasi-Pl\"ucker coordinate origin. Section \ref{sec:Qproof} shows that (\ref{eq:qComm}) does too. The second proof of Theorem \ref{thm:qComm} shows that a great many instances of $(\mathcal{M}_{J,I})$ do as well; to see this, note that the roles of $M_{J,I}$ and $C_{J,I}$ were interchangeable there. The question of whether and to what extent the gap (case $J\not\surr I$) may be filled by finding new quasi-Pl\"ucker coordinate identities is an interesting one. Toward a partial answer, we leave the reader to verify that 
\[
(\mathcal P_{I,J^j,j}) \Rightarrow (\mathcal{M}_{J,I})
\]
whenever $I,J\subseteq[n]$ are such that $|J|\leq|I|$ and $J\setminus j\subseteq I$.

\def\cprime{$'$}

\bigskip

\par\noindent
{\small\textsc{LaCIM, UQ\`AM, 
Case Postale 8888, succursale Centre-ville, \\
Montr\'eal (Qu\'ebec)\, H3C 3P8, CANADA}.\, 
\texttt{lauve@lacim.uqam.ca}}

\begin{thebibliography}{1}

\bibitem{GelRet:1}
I.~M. Gel{\cprime}fand and V.~S. Retakh.
\newblock Determinants of matrices over noncommutative rings.
\newblock {\em Funktsional. Anal. i Prilozhen.}, 25(2):13--25, 96, 1991.

\bibitem{GelRet:3}
I.~M. Gelfand and V.~S. Retakh.
\newblock Quasideterminants, {I}.
\newblock {\em Selecta Math. (N.S.)}, 3(4):517--546, 1997.

\bibitem{GGRW:1}
Israel Gelfand, Sergei Gelfand, Vladimir Retakh, and Robert~Lee Wilson.
\newblock Quasideterminants.
\newblock {\em Adv. in Math.}, 193(1):56--141, 2005.

\bibitem{KelLenRig:1}
A.~C. Kelly, T.~H. Lenagan, and L.~Rigal.
\newblock Ring theoretic properties of quantum {G}rassmannians.
\newblock {\em J. Algebra Appl.}, 3(1):9--30, 2004.

\bibitem{KroLec:1}
Daniel Krob and Bernard Leclerc.
\newblock Minor identities for quasi-determinants and quantum determinants.
\newblock {\em Comm. Math. Phys.}, 169(1):1--23, 1995.

\bibitem{Lau:1}
Aaron Lauve.
\newblock Quantum- and quasi-{P}l\"ucker coordinates.
\newblock {\em J. Algebra}, 296(2):440--461.

\bibitem{LecZel:1}
Bernard Leclerc and Andrei Zelevinsky.
\newblock Quasicommuting families of quantum {P}l\"ucker coordinates.
\newblock In {\em Kirillov's seminar on representation theory}, volume 181 of
  {\em Amer. Math. Soc. Transl. Ser. 2}, pages 85--108. Amer. Math. Soc.,
  Providence, RI, 1998.

\bibitem{TafTow:1}
Earl Taft and Jacob Towber.
\newblock Quantum deformation of flag schemes and {G}rassmann schemes, {I}. {A}
  $q$-deformation of the shape-algebra for $\mathrm{{G}{L}}(n)$.
\newblock {\em J. Algebra}, 142(1):1--36, 1991.

\end{thebibliography}
\end{document}